\documentclass[11pt]{article}   
\usepackage{amssymb,amscd,latexsym}   
\usepackage{amsmath}
\usepackage{amsthm}
\usepackage[all]{xy}
\newtheorem{Theorem}{Theorem}[section]
\newtheorem{Lemma}[Theorem]{Lemma}

\newtheorem{Proposition}[Theorem]{Proposition}

\newtheorem{Example}[Theorem]{Example}

\newtheorem{Definition}[Theorem]{Definition}

\renewcommand{\phi}{\varphi}
\newcommand{\demo}{{\sc Proof. }}

\def\qed{\hspace*{\fill} $\square$}

\def\ii{\'{\i}}
\newcommand{\xx}{{\bf x}}

\newcommand{\ff}{{\bf f}}
\renewcommand{\gg}{{\bf g}}

\newcommand{\restr}{{\kern-1pt\restriction\kern-1pt}}

\newcommand{\pp}{{\mathbb P}}

\newcommand{\cala}{{\mathcal{A}}}

\newcommand{\calo}{{\mathcal O}}

\newcommand{\cre}{{\rm Cr}}
\newcommand{\crepla}{{\rm Cr}(\PP^2)}
\newcommand{\cren}{{\rm Cr}(\PP^n)}
\newcommand{\crenn}{{\rm Cr}(\PP^{n-1})}

\newcommand{\aut}{{\rm Aut}}

\renewcommand{\star}{{\rm J}}

\newcommand{\pgl}{{\rm PGL}}

\newcommand{\jac}{{\rm Jac}}

\newcommand{\plan}{{\mathbb P}^2}

\newcommand{\PP}{{\mathbb P}}
\renewcommand{\P}{{\mathbb P}}

\newcommand{\rar}{\rightarrow}
\newcommand{\lar}{\longrightarrow}
\newcommand{\llar}{-\kern-5pt-\kern-5pt\longrightarrow}
\newcommand{\surjects}{\twoheadrightarrow}
\newcommand{\injects}{\hookrightarrow}
\newcommand{\tor}{\xymatrix{\ar@{-->}[r]&}}
\renewcommand{\tor}{\dasharrow}

\newcommand{\be}{\begin{enumerate}}
\newcommand{\ee}{\end{enumerate}}
\newcommand{\bi}{\begin{itemize}}
\newcommand{\ei}{\end{itemize}}

\newcommand{\beqn}{\begin{eqnarray*}}
\newcommand{\eeqn}{\end{eqnarray*}}

\begin{document}

\begin{center}
{\Large{\sc\sc Cremona maps of de Jonqui\`eres type}}\\
\footnotetext{2010 AMS {\it Mathematics Subject
Classification}: 13D02, 13H10, 14E05, 14E07, 14M05, 14M25.}

\vspace{0.3in}

{\large\sc Ivan Pan}\footnote{Partially supported by the \emph{Agencia Nacional de Investigaci\'on e Innovaci\'on} of Uruguay}
\quad\quad
{\large\sc Aron  Simis}\footnote{Partially
supported by a CNPq grant and by a PVNS (CAPES) Fellowship.}

\end{center}


\begin{abstract}
This paper is concerned with suitable generalizations of a plane de Jonqui\`eres map to higher dimensional space
$\pp^n$ with $n\geq 3$.
For each given point of $\pp^n$ there is a subgroup of the entire Cremona group of dimension $n$
consisting of such maps.
One studies both geometric and group-theoretical properties of this notion.
In the case where $n=3$ one describes an explicit set of generators of the group and gives a homological characterization
of a basic subgroup thereof.
\end{abstract}

\section*{Introduction}

Let $k$ denote an algebraically closed  field of characteristic zero and let $\pp^n=\pp^n_k$ denote
the projective space of dimension $n$ over $k$.
A  classical problem is to understand the structure of the $k$-automorphism group of the function field of  $\pp^n$
or, equivalently, its {\em Cremona group}  $\cre(\PP^n)$ of birational maps.

An important subgroup of $\cre(\PP^n)$ is the group $\pgl(n+1,k)$ of projective (linear) transformations.
For $n=1$ one easily sees that $\cre(\pp^1)=\pgl(2,k)$ -- the so-called group of M\"obius transformations over $k$.
For $n=2$ a celebrated result states that $\cre(\pp^2)$ is generated by $\pgl(3,k)$ and the standard quadratic map of $\PP^2.$
A first proof was given by
Max Noether (\cite{Noe71}, \cite{Noe72}).
The proof however contained one gap later filled in by Guido Castelnuovo.
For an interesting account of the history of this result, including the contributions by Castelnuovo and others, the reader
is referred to \cite[Chap. 8]{Alb}.   Now, in his version of the theorem,  Castelnuovo gives an alternative approach
by first proving that every plane Cremona transformation is a composite of  de Jonqui\`eres maps, then shows that any such map
is a composite of projective transformations and the standard quadratic map.
This shows the prominence of de Jonqui\`eres maps in the classical Cremona map theory.

A {\em de Jonqui\`eres map}
(in honor of \cite{dJo} where it has been first studied) is a plane Cremona map $\mathfrak{F}$, say, of degree $d\geq 2$,
satisfying any one of the following equivalent conditions:

\medskip

$\bullet$ $\mathfrak{F}$ has homaloidal type $(d\,; d-1,1^{2d-2})$

\smallskip

$\bullet$  There exists a point $o\in\PP^2$ such that the restriction of $\mathfrak{F}$ to a general line
passing through $o$ maps it birationally to a line passing through $o$

\smallskip

$\bullet$  Up to projective coordinate change (source and target) $\mathfrak{F}$ is defined by $d$-forms
$\{ qx_0, qx_1, f\}$ such that $f,q\in k[x_0,x_1,x_2]$ are relatively prime $x_2$-monoids one of which at
least has degree $1$ in $x_2$.

\smallskip

The first alternative means that the base cluster of the map has one proper point of multiplicity $d-1$ and $2d-2$
not necessarily distinct, possibly infinitely near points.
The second alternative is a ``dynamical'' notion emphasizing the behavior of the map with respect to proper linear subspaces
of $\pp^2$.
Finally, the third alternative stresses the shape of the defining forms or, as one might say, the underlying defining
locus of $F$, with an emphasis on the monoid shape of intervening forms.
 We refer to \cite[Proposition 2.2 and Corollary 2.3]{PanStellar} for a simple geometric proof, and to
 \cite[Proposition 2.3, also Remark 2.4]{HS} for a  later
argument stressing the algebraic fundamentals of plane Cremona maps.

Up to a projective change of coordinates, one can take $o=(0:0:1)$, and hence the explicit format delivered by the third alternative is
\[\bigl((c(x_0,x_1)x_2+d(x_0,x_1))\,x_0:(c(x_0,x_1)x_2+d(x_0,x_1))\,x_1:a(x_0,x_1)x_2+b(x_0,x_1)\bigr),\]
 $a,b,c,d\in k[x_0,x_1]$ are forms such that $ad-bc\neq 0$
and of degrees
\[deg(a)=\deg(d)=\deg(b)-1=\deg(c)+1.\]
Note that the fraction  $(a(x_0,x_1)x_2+b(x_0,x_1))/(c(x_0,x_1)x_2+d(x_0,x_1))$ defines a M\"obius transformation in the variable
$x_2$ over the function field $k(x_0,x_1)$ of $\pp^1$.

Now, in general the problems concerning the Cremona group have (at least) two facets: the group theoretic questions, such as booking
generators and relations of some important subgroups, and the geometric questions that deal with classifying types according
to the geometric properties or the group constituents.
These two facets intersperse each other, often happening that the geometric results help visualizing the group structure.
In this paper we deal with both aspects of the theory and, in addition, bring up the its underlying commutative algebra
in terms of the ideal theoretic and homological side
of the so-called base ideals of the maps.

We will deal with suitable generalizations of de Jonqui\`eres maps to higher dimensional space
$\pp^n$ with $n\geq 3$.
These generalizations will be subsumed under the general frame of {\em maps of de Jonqui\`eres type}.
For $n\geq 3$ we will study elements of the Cremona group ${\rm Cr}(n)={\rm Cr}(\pp^n)$ satisfying a condition akin to the one
of (b) above. Fixing the point $o$, these maps will form the de Jonqui\`eres subgroup ${\rm J}_o(m\,;\pp^n)\subset {\rm Cr}(n)$
of type $m$ with center $o$.
The plan is therefore to look at these two objects: the subgroup and a typical element of this subgroup.

Let us focus on the case $m=1$. The strategy for the group study is based on an exact sequence of groups
\begin{equation}\label{eq-seq}
\xymatrix{1\ar@{->}[r]&\pgl(2,k(\PP^{n-1}))\ar@{->}[r]&\star_o(1;\PP^n)\ar@{->}[r]&\cre(\PP^{n-1})\ar@{->}[r]&1},
\end{equation}
where $\pgl(2,k(\PP^{n-1}))$ corresponds to the subgroup consisting of de Jonqui\`eres transformations which stabilize
a general line passing through $o$.
The elements of this subgroup can be viewed as M\"obius transformations over the function field
$k(\PP^{n-1})$ of $\pp^{n-1}$.
Moreover, one can see that the sequence is right split (see \S \ref{sec3}).

One does not know much about the structure of $\star_o(1;\PP^n)$ when $n\geq 3$. For example, it is not known
whether this group together
with projective linear transformations generate the entire Cremona group $\cre(\PP^n)$, as it happens when $n=2$.
It has been proved in \cite[Theorem 1]{Pa99b} that if $k=\mathbb{C}$
(or, more generally, if $k$ has uncountably many elements) then any
set of group generators of $\cre(\PP^n)$  contains uncountably many non linear transformations.
In this paper we show that  $\star_o(1;\pp^3)$ itself inherits this property, but we are still
 able to describe a complete set of families of maps that generates it.

The paper is divided in three sections
The first section is devoted to the basic definitions and the main properties of the
de Jonqui\`eres subgroup ${\rm J}_o(m\,;\pp^n)\subset {\rm Cr}(n)$
of type $m$ with center $o$.

The second section establishes the main group theoretic results of the paper, based
on information coming from the geometric side.
Here we focus on dimension $3$ aiming to describe the structure of the group $\star_o(1;\pp^3)$.
The first result states that if one fixes an integer $d\geq 2$ then any $\mathfrak{F}\in \star_o(1;\PP^n)$
of degree $\leq d$ contracts a finite number of irreducible hypersurfaces each of which has geometric
genus bounded by a number depending on $\deg(\mathfrak{F})$.
This is then used to deduce that, if $k=\mathbb{C}$, given any set $\mathcal{G}$ of generators of $\star_o(1;\PP^n)$  and
an integer $d\geq 2$,
then any subset $\mathcal{G}_0\subset \mathcal{G}$ such that $\pgl(2,k(\PP^{n-1}))$ is generated by
elements of $\mathcal{G}_0$ contains uncountably many elements of degrees $\geq d$.
Finally one shows that $\star_o(1;\PP^3)$ is generated by its subgroup $\star_o(2;\PP^3)$  and by the cubic Cremona map
$(x_0x_1x_2:x_0^2x_2:x_0^2x_1:x_1x_2x_3)$.

In the third section of the paper we expand on additional algebraic aspects of rational maps akin to
maps of de Jonqui\`eres type, by stressing homological properties of a class of homogeneous ideals
resembling the base ideals of such rational maps.
Basically, bringing up an essential feature of the subject matter,  one creates a rational map
 $\pp^n\dasharrow \pp^n$ out of a given rational map  $\pp^{n-1}\dasharrow \pp^{n-1}$ and proves that the property
 of being Cremona carries over in both directions.
 In this sense, there are several ways to go from Cremona maps on $\pp^{n-1}$ to Cremona maps on $\pp^n$.
 The one studied in this paper is special in that the resulting rational map
 on $\pp^n$  automatically has base locus of codimension $2$.
Besides, one moves both source and target one dimension up and also makes the degree
of the defining forms go up as well.
A close version has been studied in \cite{HSstellar} as parametrization of certain implicit monoid hypersurfaces.
There, one moves up target by one dimension, keeping the same source.
An immediate bonus of the latter is that Cremona carries over for free to a birational map
onto its image.
However, the resulting map still has base locus of codimension $2$.
Yet a third situation of the present considerations has been dealt with in \cite{NewCremona}, where one of the goals
was to allow for a recurrently growing codimension of the base locus at the expense of not moving up
the degree of the defining forms.

The gist of the section is to take a more abstract view of the construction, with an emphasis on the ideal theoretic
and homological properties of the base locus of a rational map.
For most of the material of this part one can drop the requirement on the characteristic of $k$
and any additional hypothesis on the transcendence degree of $k$ over its prime field.

\section{Maps of the Jonqui\`eres type}

We will be solely concerned with rational maps of $\pp^n$ to itself.
A rational map $\mathfrak{F}:\pp^n\dasharrow \pp^n$ is defined by $n+1$ forms $\ff=\{f_0,\ldots, f_n\}
\subset R:=k[\xx]=k[x_0,\ldots,x_n]$ of the same degree $d\geq 1$, not all null.
We often write $\mathfrak{F}=(f_0:\cdots :f_n)$ to underscore the projective setup.
Any rational map can without lost of generality be brought to satisfy  the condition
that $\gcd\{f_0,\cdots ,f_n\}=1$ (in the geometric terminology, $\mathfrak{F}$ {\em has no fixed part}).
In order to have a well-defined notion of {\em degree} of $\mathfrak{F}$ we will always assume the latter condition
which means that we will usually be identifying the maps $(f_0:\cdots :f_n)$ and $(ff_0:\cdots :ff_n)$
for any nonzero form $f\in R$.

Given a point $o\in\PP^n$, we denote by $S_o(m,n)$ the Schubert cycle parameterizing the set of
$m$-dimensional linear subspaces
$L\subset\PP^n$ containing $o$; it is known to be an irreducible variety.
We will as usual identify a member $L$ of the set with the corresponding point of $S_o(m,n)$.

Fix a point $o\in\PP^n$ and an integer $1\leq m\leq n-1$.
\begin{Definition}\rm A Cremona map of $\pp^n$ is a {\em  de Jonqui\`eres map of type $m$ with center} $o$ if there is a
sufficiently small nonempty open set $U\subset \pp^n$ such that the restriction of $F$ to $U$ is biregular onto its image
and for which the following condition holds: given $L\in S_o(m,n)$ such that $L\cap U\neq\emptyset$ then
$\overline{F(L\cap U)}\in S_o(m,n)$, where over-line indicates Zariski closure.
\end{Definition}
In a more informal way, the condition is that the general member of the the Schubert cycle $S_o(m,n)$ is mapped onto a
member of $S_o(m,n)$. Note that since the set of contracted linear subspaces by a Cremona map is contained in a finite set of hypersurfaces
then the restriction of the Cremona map to a general member of $S_o(m,n)$ is a birational
map onto its image. Moreover, since a general member of $S_o(m,n)$ is the intersection of two general members in
$S_o(m+1,n)$ one easily deduces that ${\rm J}_{o}(m'\,;\pp^n)\subset {\rm J}_{o}(m\,;\pp^n)$
for any $m'\geq m$.

\bigskip

The set ${\rm J}_o(m\,;\pp^n)$ of all de Jonqui\`eres maps of type $m$ with center $o$ is a subgroup
of the whole Cremona group -- it will be referred to as the {\em de Jonqui\`eres group of type $m$ with center} $o$.
We note that, with varying terminology, this notion has appeared elsewhere before -- e.g., \cite[Proposition 2.1]{PanStellar}, also
\cite[Section 4.3]{Dol} where the author called {\em level} $n-m$ what we call {\em type} $m$.

We focus on the case $m=1$.
Let $H\subset\PP^n$ be a hyperplane not containing the point $o$ and let $K$ stand for its field of rational functions over $k$; note that the projective space of lines passing through $o$ may be identified with $H$ by associating each such line to its intersection with $H$. 
By definition, an element $\mathfrak{F}\in\star_o(1;\PP^n)$ is a Cremona map of $\PP^n$ that acts birationally
on the set of lines passing through $o$, hence $\mathfrak{F}$ induces a  birational maps $H\dasharrow H$. 
By identifying $H=\PP^{n-1}$, we obtain a map
\begin{equation}\label{main_map}
\rho:\star_o(1;\PP^n)\to \cre(\PP^{n-1})
\end{equation}
which is clearly a group homomorphism.

The group $\star_o(1;\PP^n)$ itself has been  treated in \cite{PanStellar}.
Here we provide further details about the above map.

For this, we introduce some additional notation, where we set $o=(0:\cdots :0:1)$
and $H:\{x_n=0\}$.
Given forms $a,b,c,d\in k[x_0,\ldots, x_{n-1}]={\rm Sym}_k[H^*]$ of degrees $r-1, r, r-2, r-1\geq 1$, respectively,
such that either $a\neq 0$ or $c\neq 0$ and satisfying $\gcd(ax_n+b,cx_n+d)=1$, one consideres the following
two objects:

\begin{itemize}
\item The element $\mathfrak{f}_{a,b,c,d}:=(ax_n+b)/(cx_n+d)\in k[x_0,\ldots, x_{n-1}](x_n)$
\item The rational map $\mathfrak{F}_{a,b,c,d}:\pp^n\dasharrow \pp^n$ {\rm (}$n\geq 2${\rm )}  defined by
\begin{equation}\label{proped}
\left( (cx_n+d)x_0:\cdots : (cx_{n}+d)x_{n-1}: ax_n+b\right)
\end{equation}
\end{itemize}
Let us emphasize that always $ax_n+b\in (x_0,\ldots, x_{n-1})k[x_0,\ldots, x_{n-1}][x_n]$, while
the same holds for $cx_n+d$ except when its degree as a form is $1$.
We observe that $\mathfrak{F}_{a,b,c,d}$ is an element of $\cren$ since the first $n$ coordinates
defines the identity map of $\cre(H)$ up to the identification $H={\rm Proj}({\rm Sym}_k(k[H^*]))$
and the last coordinate is of degree $\leq 1$ in the
variable $x_n$ (\cite[Proposition 2.2]{PanStellar}).

Forms such as $ax_n+b, cx_n+d$ are called $x_n$-{\em monoids}.

\begin{Proposition}\label{cremona_one_variable}
Let $\rho:\star_o(1;\PP^n)\to \cre(\PP^{n-1})$ be as in {\rm (\ref{main_map})}.
Fix $o=(0:\cdots :0:1)$ and $H:\{x_n=0\}$ and let $K$ stand for the function field
of $H\subset \pp^n$.
Then
\begin{enumerate}
\item[{\rm (i)}] A Cremona map $\mathfrak{F}\in\cren$ belongs to $\star_o(1;\PP^n)$ if and only if
as a rational map $\PP^n \dasharrow \PP^n$ it has the form
\begin{equation}
(qg_0:\cdots:qg_{n-1}:f),
\end{equation}
where $(g_0:\cdots:g_{n-1})$
defines a Cremona map of $H\simeq \PP^{n-1}$ and $q,f\in k[x_0,\dots,x_{n-1},x_n]$ are
relatively prime $x_n$-monoids one of which at least has positive $x_n$-degree.
\item[{\rm (ii)}] The group $\pgl\left(2,K\right)$ can be identified with the M\"obius group whose
elements have the form $\mathfrak{f}_{a,b,c,d}$.
\item[{\rm (iii)}] The map $\mathfrak{f}_{a,b,c,d}\mapsto \mathfrak{F}_{a,b,c,d}$ is an injective
group homomorphism $\psi: \pgl\left(2,K\right) \injects \star_o(1;\PP^n)$.
\item[{\rm (iv)}] ${\rm im}(\psi)=\ker(\rho)\,${\rm ;} in particular $\mathfrak{F}_{a,b,c,d}$ maps a
general line passing through the point $o=(0:\cdots:0:1)$ birationally to itself.
\end{enumerate}
\end{Proposition}
\demo
(i) This has been proved in \cite[Proposition 2.2]{PanStellar}, but we repeat the main steps for the
reader convenience. (A different proof stressing the field theoretic side is available but will not be 
given to avoid distraction.) 
Consider the coordinate  projection  $\pi:\PP^n\tor H=\{x_n=0\}$  with center $o=(0:\dots:0:1)$:
\[(x_0:\dots:x_n)\mapsto (x_0:\dots:x_{n-1}).\]
 Then a Cremona map $\mathfrak{F}\in\cren$
belongs to $\star_o(1;\PP^n)$ if and only if there exists $\mathfrak{G}\in\cre(H)=\cre(\PP^{n-1})$ such that
$\pi\mathfrak{F}=\mathfrak{G}\pi$.
It follows that $\mathfrak{F}\in\star_o(1;\PP^n)$ if and only if  there exist forms $q,f\in k[x_0,\dots,x_n]$ of
respective degrees $\deg(\mathfrak{F})-\deg(\mathfrak{G})$ and  $\deg(\mathfrak{F})$, such that
\[\mathfrak{F}=(qg_0:\cdots:qg_{n-1}:f)\]
for a suitable  $\mathfrak{G}=(g_0:\cdots:g_{n-1})\in\cre(H)\simeq \cre(\PP^{n-1})$,
where $g_i\in k[x_0,\dots ,x_{n-1}]$ for $i=0,\dots,n-1$.
To see that $q,f$ are necessarily $x_n$-monoids, one uses the implication (1) $\Rightarrow$ (2) of
\cite[Proposition 2.2]{PanStellar}.

(ii) This is an easy exercise passing to inhomogeneous coordinates
$\frac{x_1}{x_0},\ldots,\frac{x_{n-1}}{x_0}.$

(iii)
By (i), $\psi$ maps to $\star_o(1;\PP^n)$.
A straightforward computation gives the composition law
$$\mathfrak{F}_{a,b,c,d}\,\mathfrak{F}_{a',b',c',d'}=\mathfrak{F}_{aa'+bc',ab'+bd',ca'+dc',cb'+dc'},$$
which shows that $\psi$ is a group homomorphism.
For the  injectivity note that  $\mathfrak{F}_{a,b,c,d}$ is the identity map if and only if $b=c=0, a=d$.

(iv) Clearly $\rho$ maps any $\mathfrak{F}_{a,b,c,d}$ to the identity map of $\cre(\PP^{n-1})$.
Conversely, let $\mathfrak{F}\in\ker(\rho)$.  By (i),  $\mathfrak F=(qg_0:\cdots:qg_{n-1}:f)$,
for suitable $x_n$-monoids $q,f\in k[x_0,\dots,x_{n-1},x_n]$ one of which at least has positive $x_n$-degree.
But since $\rho$ maps $\mathfrak{F}$ to the identity of $\cre(\PP^{n-1})$ then $(t_0:\cdots t_{n-1})$ must
be the identity map.
This shows that $\mathfrak{F}=\mathfrak{F}_{a,b,c,d}$, with $f=ax_n+b, q=cx_n+d$.
\qed

\medskip

\medskip

Let  $o=(0:\dots:0:1)$ and $H=\{x_n=0\}\subset \pp^n$ as before.

As a consequence of the above methods, we observe that  $\star_o(1;\PP^n)$
has two distinguished subgroups: one is the kernel $\ker(\rho)$
which we have shown to be exactly the subgroup of Cremona maps of the form $\mathfrak{F}_{a,b,c,d}$, for suitable forms
$a,b,c,d\in k[x_0,\ldots,x_{n-1}]$.
Note that these fix a general hyperplan through the point $o$ since the first $n$ coordinates of the map
define the identity map on the fixed hyperplan $H$ avoiding $o$.
The other  subgroup is $\star_o(n-1;\PP^n)\subset \star_o(1;\PP^n)$ whose elements
map a general hyperplan through $o$ birationally onto a hyperplan through $o$ (not necessarily
fixing the source hyperplan).
The simple geometry behind the relationship of these two subgroups asks for a group-theoretic
formulation.
And in fact, there is a simple one:

\begin{Proposition}\label{pro_fact}
Let $\pgl(n+1,k)_o$ be the subgroup of linear automorphisms of $\PP^n$  fixing $o$. Then
$\ker(\rho)$ is a normal subgroup of $\star_o(n-1;\PP^n)$ and the equality
\[\star_o(n-1;\PP^n)=\pgl(n+1,k)_o\ker(\rho)\]
holds.
\end{Proposition}
\demo
Clearly, both $\ker(\rho)$ and $\pgl(n+1,k)_o$ are subgroups of $\star_o(n-1;\PP^n)$ and the first
is normal since it is normal in the larger group $\star_o(1;\PP^n)$. Therefore the product of the two subgroups
is a subgroup of $\star_o(n-1;\PP^n)$.
Conversely, let $\mathfrak{F}\in \star_o(n-1;\PP^n)$. Expressing it as an element of $\star_o(1;\PP^n)$
we know from the previous part that $\mathfrak{F}=(qg_0:\dots:qg_{n-1}:f)$,
for a suitable  $\mathfrak{G}=(g_0:\cdots:g_{n-1})\in\cre(\PP^{n-1})$ and certain $x_n$-monoids $q,f$.
 Since $\mathfrak{F}$ maps a general hyperplane through $o$ to a hyperplane
through $o$, the forms $g_0:\cdots:g_{n-1}$ are necessarily linear forms in $k[x_0,\ldots,x_{n-1}]$.
Let $ \mathfrak{A}\in\pgl(n+1,k)_o$ denote the inverse of the linear automorphism defined by
$(g_0:\cdots :g_{n-1}:x_n)$.
Then $\mathfrak{A}\mathfrak{F}=(qx_0:\dots:qx_{n-1},f)$, which is a map of the type
$\mathfrak{F}_{a,b,c,d}$, and hence it belongs to $\ker(\rho)$ by
Proposition~\ref{cremona_one_variable} (iv).
\qed

\section{Generators of the de Jonqui\`eres group of type $1$}\label{sec3}

The following fact has been established in \cite[Proposition 2.1]{PanStellar}, for which an affine argument was given.
We isolate it as a lemma for reference convenience and give a proof in terms of the projective geometry.

\begin{Lemma}\label{semi_direct}
Consider the previous group homomorphism 
$\star_o(1;\PP^n)\stackrel{\rho}{\lar}\cre(\PP^{n-1}),$
whose kernel is identified with  $\pgl(2,k(\PP^{n-1}))$ by 
{\rm Proposition~\ref{cremona_one_variable}}.
Then there is a map $\sigma:\crenn\to \star_o(1;\PP^n)$ such that $\rho\circ\sigma$ is the
identity of $\cre(\PP^{n-1})$.
In particular, $\rho$ is surjective ans $\star_o(1;\PP^n)$ is isomorphic to the semi-direct product $\pgl(2,k(\PP^{n-1}))\rtimes \crenn$.
\end{Lemma}
\demo A splitting map $\sigma$ is of course not uniquely defined.
We choose one such map $\sigma:\crenn\to \star_o(1;\PP^n)$: given  $\mathfrak{t}:=(t_0:\dots:t_{n-1})\in\crenn$,
let $\sigma(\mathfrak{t})$
be the rational map of $\pp^n$ defined as follows:
\begin{equation}\label{eq-semi}
\sigma(\mathfrak{t})=(x_0t_0({\mathbf x}):\dots:x_0t_{n-1}({\mathbf x}):t_0({\mathbf x})\,x_n),
\end{equation}
where ${\mathbf x}=\{x_0,\dots ,x_{n-1}\}$.
It is clear that $\sigma(\mathfrak{t})\in \star_o(1;\PP^n)$ by appealing to Proposition~\ref{cremona_one_variable} (i)
with $\mathfrak{G}=(t_0({\mathbf x}):\dots:t_{n-1}({\mathbf x}))$, $q=x_0$ and $f=t_0({\mathbf x})\,x_n$.

Note that if $\mathfrak{t}'=(t'_0:\dots:t'_{n-1})\in\crenn$ has degree $r$ and $\mathfrak{t}\circ\mathfrak{t}'=(s_0:\dots:s_{n-1})$, then
\begin{eqnarray*}
\sigma(\mathfrak{t})\circ\sigma(\mathfrak{t}')&=&\sigma(\mathfrak{t})\bigl(x_0t'_0({\mathbf x}),\dots,x_0t'_{n-1}({\mathbf x}),x_nt'_0({\mathbf x})\bigr)\\
&=&(x_0^{r+1}s_0({\mathbf x}):\dots:x_0^{r+1}s_{n-1}({\mathbf x}):x_0^rx_ns_0({\mathbf x}))\\
&=&(x_0s_0({\mathbf x}):\dots:x_0s_{n-1}({\mathbf x}):x_ns_0({\mathbf x}))\\
&=&\sigma(\mathfrak{t}\circ\mathfrak{t}').
\end{eqnarray*}
By definition of the map $\rho:\star_o(1;\PP^n)\to\cre(\PP^{n-1})$ it is clear that $\rho\circ\sigma$ is the
identity map of the group $\cre(\PP^{n-1})$.
\qed

\medskip

We proceed to elucidate further the geometric behavior.
For this recall that a rational map $\mathfrak{F}: \pp^n\dasharrow \pp^n$ {\em contracts} a subvariety $V\subset \pp^n$
provided the restriction of $\mathfrak{F}$ to an open dense subset of $V$ is well-defined and its image has dimension
strictly less than $\dim V$.
It is known that if $\mathfrak{F}$ contracts an irreducible hypersurface $V\subset \pp^n$ then the defining
equation of $V$ is a factor of the Jacobian determinant of the forms defining $\mathfrak{F}$.

If $V\subset\PP^r$ is an irreducible  projective variety of dimension $m$, the \emph{geometric genus} $p_g(V)$ of $V$ 
is the maximal number of linearly independent global differential $m$-forms on some (then all)  desingularization of $V$;
if $V$ is a smooth hypersurface,  then $p_g(V)={\ell -1\choose r-1}$. In general, writing $\ell-1=s(r-m)+e$, where $0\leq e\leq r-m-1$, one has the 
so-called Castelnuovo--Harris bound for the geometric genus of $V$ (see \cite{Ha}):
\[p_g(V)\leq {s \choose m+1}(r-m)+{s\choose m}e.\]

 \begin{Proposition}\label{contraction}
Let $d\geq 2$ be a given integer.
Then any $\mathfrak{F}\in \star_o(1;\PP^n)$ of degree $\leq d$ contracts
a finite number of irreducible hypersurfaces each of which has geometric genus bounded by a number
depending on $\deg(\mathfrak{F})$.
\end{Proposition}
\demo
Let $\mathfrak{F}\in\star_o(1;\PP^n)$ be a de Jonqui\`eres map of type $1$.
By Proposition~\ref{cremona_one_variable} (i)
there exists $\mathfrak{G}=(g_0:\dots:g_{n-1})\in\crenn$ such that
\begin{equation}\label{format_repeated}
\mathfrak{F}=(qg_0:\dots:qg_{n-1}:f),
\end{equation}
for $x_n$-monoids  $q=cx_n+d, f=ax_n+b$, for suitable forms $a,b,c,d\in k[x_0,\ldots,x_{n-1}]$, with either $c\neq 0$ or $a\neq 0$ and
and $\gcd(q,f)=1$ (equivalently, $ad-bc\neq 0$); set $\deg \mathfrak{G}=\deg \mathfrak{F}-t$, 
where $t=\deg(q)$ with $1\leq t\leq \deg \mathfrak{F}$.

Note that $\mathfrak{F}$ maps  the hyperplane $x_n=0$ birationally onto a hypersurface. We deduce that $\mathfrak{F}$ is a 
local isomorphism at a point $p\in\PP^n$ if and only if $p$ is not a zero of
$(ad-x_n^{\deg\mathfrak{F}-t}bc)q\jac(\mathfrak{G})$, where $\jac(\mathfrak{G})$ is the Jacobian determinant of the set 
$\{g_0,\dots,g_{n-1}\}$: indeed, for $p$ belonging to the open set $\{x_n=1\}$, if $q(p)\neq 0$, then that map is not a 
local isomorphism at $p$ if and only if it is a zero of the
Jacobian determinant of the set $\{g_0,\dots,g_{n-1},f/q)\},$
i.e., a zero of $(ad-bc)\jac(\mathfrak{G})$. Then 
the reduced hypersurface $J(\mathfrak{F})$ of equation  $(ad-x_n^{\deg\mathfrak{F}-t}bc)q\jac(\mathfrak{G})=0$
 has degree $\leq t^2(\deg\mathfrak{F}-1)n(\deg \mathfrak{G}-1)\leq n(\deg \mathfrak{F})^2(\deg \mathfrak{F}-1)^2$.

Now let $V\subset\PP^n$ denote an irreducible hypersurface contracted by $\mathfrak{F}$.
One knows that $V\subset J(\mathfrak{F})$.  Then the assertion follows by using the Castelnuovo--Harris bound for $V$.
\qed

\medskip

For the next result we assume that the ground field is uncountable (e.g., the field of complex numbers).

\begin{Theorem}\label{thm1}
Let $\mathcal{G}$ be a set of generators for $\star_o(1;\PP^n)$  and let $d\geq 2$.
Then any subset $\mathcal{G}_0\subset \mathcal{G}$ such that $\pgl(2,k(\PP^{n-1}))$ is generated by
elements of $\mathcal{G}_0$ contains uncountably many elements of degrees $\geq d$.
\end{Theorem}
\demo
First note that for every $\ell$, and every irreducible smooth hypersurface
$\Gamma\subset \PP^{n-1}$ of degree $\ell$, we may construct an element of $\ker(\rho)\subset \star_o(1;\PP^n)$
of degree $\ell+1$ following the recipe in (\ref{format_repeated}):  take $\mathfrak{G}$ to be the identity map,
$q\in k[x_0,\dots,x_{n-1}]$ to be a defining equation of $\Gamma$ and $f=ax_n+b$, with $a\neq 0$ .

Recall  that if $\ell>n$, then  two smooth hypersurfaces $\Gamma_1, \Gamma_2\subset \PP^{n-1}$ of degree $\ell$,
which are birationally equivalent, are necessarily biregularly so --  this is because the canonical class of such a
hypersurface is ample  (cf. \cite[Thm. 0.2.1]{Che}). In particular, if $\Gamma_1$ and $\Gamma_2$ are non isomorphic, then $\Gamma_1\times\PP^1$ and $\Gamma_2\times\PP^1$ are no birationally equivalent by L\"uroth's Theorem. 

Since the moduli space  of smooth hypersurfaces of degree $\ell$ in $\PP^{n-1}$ has positive dimension, then there are 
uncountably many such hypersurfaces which are pairwise non isomorphic.  We deduce that for any $\ell>n$ there exists 
a family $\cala_\ell\subset \star_o(1,\PP^n)$ with the following properties:

\begin{itemize}

\item $\cala_\ell$ contains uncountably many elements.

\item if $\mathfrak{F}\in\cala_\ell$ then $\deg(\mathfrak{F})=\ell+1$ and $\mathfrak{F}\in\ker(\sigma)$

\item if $\mathfrak{F}_1, \mathfrak{F}_2\in\cala_\ell$, then there are irreducible components $V_1\subset J(\mathfrak{F}_1)$ and
$V_2\subset J(\mathfrak{F}_2)$, which are not birationally equivalent: indeed, we may choose $V_i\subset\PP^n$ to be a cone, with vertex $o$, over a smooth hypersurface $\Gamma_i\subset \{x_0=0\}=\PP^{n-1}$ of degree $\ell$, then $V_i$ is birationally equivalent to $\Gamma_i\times \PP^1$.

\end{itemize}

It follows from  (\cite[Lemma 4]{Pa99b}) that if $\mathcal{G}_1$ is a subset of $\cre(\PP^n)$ and if $\mathfrak{F}\in \cre(\PP^n)$
is written as a product of elements in $\mathcal{G}_1$, then every hypersurface contracted by $\mathfrak{F}$ is birationally  
equivalent to a hypersurface contracted by some element in  $\mathcal{G}_1$.

Suppose that there exist an integer $d\geq 2$ and a subset $\mathcal{G}_0\subset \mathcal{G}$ such that $\mathcal{G}_0$ contains at most 
countably many elements of degree $\geq d$. Now every member of such a countable sequence of elements of degrees $d_e:=d+e$, $e=0,1,\ldots$, contracts a 
finite number of hypersurfaces. We then deduce that, for any  $\ell > n$, an uncountable subset $\cala^{'}_\ell\subset \cala_\ell$ is generated 
by elements of degree $\leq d$ belonging to  $\mathcal{G}_0$.  This contradicts Proposition~\ref{contraction} since it suffices to take 
$\ell$ large enough in order to obtain elements in $\cala^{'}_\ell$ which contract a hypersurface with geometric genus 
larger than the stated bound.
\qed

\medskip

We now turn to the question of giving explicit families of maps that together generate $\star_o(1;\PP^n)$.
Since $\star_o(1,\PP^n)$ is isomorphic to $ \pgl(2,k(\PP^{n-1})\rtimes \crenn$, it suffices to find generators for the two factors.
In the case where $n>3$, however, this approach is helpless since we do not know a workable set of generators of $\crenn$.
On the other hand, for $n=3$ the set of generators found this way is not so exciting
since $\cre(\pp^2)$ is sufficiently familiar.
We will nevertheless state the explicit result in the next proposition.
For this, recall the generation of the M\"obius group $\pgl(2,K)$, where $K$ is an arbitrary field,
by the elements
defined by matrices of the following types
\[\begin{pmatrix}
\alpha& 0\\
0&1
\end{pmatrix}, \begin{pmatrix}
1& \beta\\
0&1
\end{pmatrix}, \begin{pmatrix}
0& 1\\
1&0
\end{pmatrix}\]
where $\alpha,\beta\in K$, $\alpha\neq 0$.
These three types of matrices are often called \emph{elementary Moebius maps} over $K$,
generating, respectively, the torus $K^*$, the additive group $K$ and the order 2 cyclic group defined by the
``inversion'' $t\mapsto 1/t$.

\begin{Proposition}\label{pro2.2}
The group $\star_o(1;\PP^3)=\pgl(2,k(\PP^2))\rtimes\crepla$ is generated by the elements $1\rtimes\mathfrak{s}$ and
$1\rtimes\mathfrak{t}$, $\mathfrak{f}\rtimes 1$, where $\mathfrak{s}$ is the standard quadratic plane Cremona map,
 $\mathfrak{t}\in\pgl(3,k)$, and $\mathfrak{f}\in \pgl(2,k(\plan))$ is an elementary Moebius map.
\end{Proposition}
\demo
The result is immediate from the above prolegomena on the generation of the M\"obius group and
the Noether--Castelnuovo Theorem on the generation of $\cre(\pp^2)$ by the standard quadratic transformation and  by $\pgl(3,k)$.
\qed

\medskip

An internal description of a set of generators of $\star_o(1;\PP^3)$ as a subgroup of $\cre(\PP^3)$
will be a consequence of Proposition~\ref{pro2.2}, through the required interpretation.
We can restate the proposition in the following compact form.

\begin{Theorem}\label{generation}
$\star_o(1;\PP^3)$ is generated by $\star_o(2;\PP^3)$  and by the cubic Cremona map
\begin{equation}\label{the_cubic}
\mathfrak{T}_3=(x_0x_1x_2:x_0^2x_2:x_0^2x_1:x_1x_2x_3)
\end{equation}
\end{Theorem}
\demo
Note that $\mathfrak{T}_3$ belongs to $\star_o(1;\PP^3)$, with $\mathfrak{G}$ the standard quadratic map
of $\cre(\pp^2)$, $q=x_0$ and $f=(x_1x_2)x_3$ (clearly, $ad-bc=x_0x_1x_2\neq 0$).
As an element of the semi-direct product decomposition it is $1\rtimes \mathfrak{s}=\sigma(\mathfrak{s})$,
clearly an involution since $\mathfrak{s}$ is.

On the other hand, any $\mathfrak{t}\in\pgl(3,k)=\aut(\PP^2)$ is defined by three independent linear forms $\ell_0,\ell_1,\ell_2\in k[x_0,x_1,x_2]$.
Its image by $\sigma$ is $1\rtimes \mathfrak{t}=(x_0\ell_0:x_0\ell_1:x_0\ell_2:x_3\ell_0)$, which has degree $\leq 2$ (actually, $=2$ provided
$\ell_0\neq \alpha x_0$, for $\alpha\in k$).
Note that $1\rtimes \mathfrak{t}$ belongs to $\star_o(2;\PP^3)$.

Finally, we know that the $\pgl(2,k(\plan))$ is identified with $\ker(\rho)\subset \star_o(2;\PP^3)$.
\qed

\medskip

To close the section, we state yet another result individualizing further the set of generators into some explicit families.

For this, we consider the set ${\bf T_{22}}$ of Cremona maps of $\pp^3$ of degree $2$ with inverse of degree $2$.
We consider the usual action of $\pgl(4,k)\times\pgl(4,k)$ on ${\bf T_{22}}$:
\[(\mathfrak{T}_1,\mathfrak{T}_2)\cdot \mathfrak{F}=\mathfrak{T}_1\mathfrak{F}\mathfrak{T}_2^{-1}.\]
One knows from \cite[Prop. 2.4.1 and Thm. 3.1.1]{PRV}  that  ${\bf T_{22}}$ is an irreducible variety of dimension $26$
with $7$ orbits under the natural action of $\pgl(4,k)\times \pgl(4,k)$.
It can be seen that any $\mathfrak{F}\in{\bf T_{22}}$ admits points
$o_1,o_2\in\PP^3$ such that $\mathfrak{F}$ transforms a general plane going through $o_1$ in a plane going through $o_2$.
Up to projective transformations  every orbit meets $\star_o(2;\PP^3)$, and also
$\ker(\rho)\subset \star_o(1;\PP^3)$ by Proposition~\ref{pro_fact}.

Considering the induced action of $\pgl(4,k)_o\times\pgl(4,k)_o$  on $\ker(\rho)\cap {\bf T_{22}}$, one has:

\begin{Lemma}
Any orbit of  $\ker(\rho)\cap {\bf T_{22}}$ under the action of $\pgl(4,k)_o\times\pgl(4,k)_o$ is the restriction to $\ker(\rho)$
of an orbit of ${\bf T_{22}}$ under the action of $\pgl(4,k)\times \pgl(4,k)$. In addition,  every such orbit of
$\ker(\rho)\cap {\bf T_{22}}$ contains involutions.
\end{Lemma}
\demo
Let $O\subset{\bf T_{22}}$ be an orbit under the action of $\pgl(4,k)\times\pgl(4,k)$. We know that $O\cap \ker(\rho)\neq\emptyset$.
To prove the first assertion it suffices to
show that for $\mathfrak{F},\mathfrak{G}\in O$ and given $\mathfrak{T},\mathfrak{T}'\in\pgl(4,k)$ such that
$\mathfrak{T}\mathfrak{F}\mathfrak{T}'=\mathfrak{G}$, then $\mathfrak{F}, \mathfrak{G}\in \ker(\rho)$ implies $\mathfrak{T},\mathfrak{T}'\in\pgl(4,k)_o$.

The linear system associated to an element in ${\bf T_{22}}$ is defined by smooth quadrics containing a conic $C$ of rank $r$,
where $r\in\{1,2,3\}$, and going through a unique ``special'' point $p$ which either does not belong to the conic plane or $p\in C$
and the tangent plane at $p$ of a general member in that linear system is constant and does not contain $C$.
In particular, such an element in ${\bf T_{22}}$ is not defined (only) along $C\cup\{p\}$ and if it belongs to $\pgl(4,k)\star_o(1;\PP^3)$ then $p=o$.

Let $L\simeq\P^1$ be a line in $\PP^3$.  If $L$ intersects the open set on which $\mathfrak{F}$ is injective, the restriction of $\mathfrak{F}\in \ker(\rho)\cap {\bf T_{22}}$ to $L$
induces a biregular map onto the image $\nu:\PP^1\to \PP^3$.  Since $\mathfrak{F}$ is defined by quadratic polynomials,
then $\nu^*\calo_{\PP^3}(1)=\calo_{\PP^1}(n)$ with  $n\in\{1,2\}$. We have $n=1$ if and only if all these polynomials
vanish at a point $p'\in L$; in this case  the point $p'$ is the special point associated to $\mathfrak{F}$ as element in ${\bf T_{22}}$, that is, $p'=o$.

By applying this argument to lines of the form $L=\mathfrak{T}'(L_o)$, where $L_o$ is a general line passing through $o$, and taking into account that $\mathfrak{G}$ maps $L_o$ birationally onto a line of the same type,
we deduce that the special point of $\mathfrak{T}\mathfrak{F}$ is $o$, hence $\mathfrak{T}'(o)=o$. By symmetry the same
holds for $\mathfrak{T}^{-1}(L_o)$ and  $(\mathfrak{T}')^{-1}\mathfrak{F}^{-1}$, hence $\mathfrak{T}(o)=o$.

This proves the first statement.
The second assertion may be proved using the normal forms obtained in \cite[Thm. 3.1.1]{PRV} or drawing upon the
content of \cite[Cor. 5.3]{PiRu} (also \cite[Thm. 5.11]{PiRu}).
\qed

\medskip

We now deduce:

\begin{Theorem}\label{thm2.3}  $\star_o(1;\PP^3)$ is generated by the cubic Cremona involution $\mathfrak{T}_3$,
seven  involutions of degree $2$,
generators for $\pgl(4,k)_o$  and the Cremona maps of degree $\geq 3$ coming from elementary Moebius map in $\pgl(2,k(\PP^2))$.
\end{Theorem}

\section{Generalized de Jonqui\`eres ideals}\label{main_de_Jonq}

In previous sections we focused on the group theoretic and geometric properties of
certain Cremona maps. As we have seen, these maps admit a very special form in terms of their defining coordinates.
In this section we take a more abstract view of a map of de Jonqui\`eres type, with an emphasis on the ideal theoretic
and homological properties of the base ideal thereof.

Quite generally, let there be given a rational map
\begin{equation}\label{base_map}
\mathfrak{G}=(g_0:\cdots :g_{n-1}):\pp^{n-1}_k\dasharrow\pp^{n-1}_k,
\end{equation}
where $g_i$'s are forms of degree $d\geq 1$ in the polynomial ring $R:=k[\xx]=k[x_0,\ldots, x_{n-1}]$.
Write $I:=(g_0,\ldots,g_{n-1})\subset R$ for the corresponding base ideal.
Consider the flat extension $S:=k[\xx,x_n]=R[x_n]=k[x_0,\ldots, x_{n-1},x_n]$, where $x_n$ is a new
indeterminate.
Let $q,f\in S$ be additional forms of degrees  $\mathfrak{d}\geq 1$ and $D:=d+ \mathfrak{d}$,
respectively, where $\mathfrak{d}$ is arbitrary.
We assume throughout that $q$ and $f$ are relatively prime.

\begin{Definition}\label{q,f-type}\rm
The rational map $\mathfrak{F}:=(qg_0:\cdots :qg_{n-1}:f):\pp^n\dasharrow\pp^{n}$ will be
called a {\em  map of $(q,f)$-type} with {\em underlying map} $\mathfrak{G}$.
\end{Definition}
The corresponding base ideal $J:=(qI,f)\subset R$ will also be called a {\em $(q,f)$-ideal}
with {\em underlying ideal} $I$.

\subsection{The homology of a map of $(q,f)$-type}\label{2.1}

As remarked above, beyond the formal similarity between the setup of \cite{PanStellar} and that of \cite{HSstellar},
the two carry substantially different contents since
the second is always  a birational map onto its image (``the implicit equation''),
while the first is not Cremona unless further constraints are put on the nature of the forms $q, f$
in the definition.
Nevertheless, the two share a common homological thread, which we would like to explore.

One reason to explore this homological facet is the following result:
\begin{Proposition}
Let $I\subset R:=k[x_0,\ldots, x_{n-1}, x_n]$
denote the base ideal of an arbitrary map in $\star_o(1;\PP^n)$ that maps a general line passing through $o$
birationally to itself.
Then $I$ is a perfect ideal {\rm (}i.e., $R/I$ is a Cohen--Macaulay ring{\rm )} if and only if $n=2$.
\end{Proposition}
\demo
By the previous proposition, we know that $\mathfrak{F}=\mathfrak{F}_{a,b,c,d}$, for suitable forms
$a,b,c,d\in k[x_0,\ldots,x_{n-1}]$.
It follows that the following matrix is part of a matrix of syzygies of the generators of $I$:
\begin{equation}\label{koszul}
\phi= \left(
\begin{array}{cc}
\kappa & \mathfrak{c}(ax_n+b)\\
\mathbf{0} & cx_n+d
\end{array}
\right).
\end{equation}
Here $\kappa$ denotes the matrix of the Koszul syzygies of the variables $x_0,\ldots, x_{n-1}$
and $\mathfrak{c}(ax_n+b)$ is the column vector whose coordinates are the content coefficients of $ax_n+b$
as an element of the extended ideal $(x_0,\ldots, x_{n-1})R$ --
for convenience the rightmost syzygy is referred to as the {\em content syzygy} of $I$
(actually, $\phi$ will be the entire syzygy matrix as will become clear in the third section).

As an immediate consequence of the size of the matrix (\ref{koszul}),  which has ${{n}\choose {2}}+1$
independent columns, if $I$ is perfect then it must be the case that${{n}\choose {2}}+1=n$.
The only possibility with $n\geq 2$ is $n=2$.
Conversely, if $n=2$ with $\mathfrak{F}$ of the form $\mathfrak{F}_{a,b,c,d}$, it is immediate to see that
the coordinates of the map are the $2$-minors of a $3\times 2$ matrix, hence generate a perfect ideal.
\qed

\medskip

Thus, one sees that there is some nontrivial homological behavior behind the scenes.

\smallskip

Recall that an $x_n$-monoid is a form $a(\xx)x_n+b(\xx)$, where $a(\xx),b(\xx)$ are forms in $k[\xx]=k[x_0,\ldots,x_{n-1}]$.
If $a(\xx)\neq 0$, such a form is irreducible if and only if $\gcd(a,b)=1$.
We keep the convention that, by definition, a Cremona map has base locus of codimension $\geq 2$.

We keep the notation of ({\ref{base_map}) and Definition~\ref{q,f-type}.

These rational maps admit a fairly structured homological behavior.

Let as before $I\subset R=k[x_0,\ldots,x_{n-1}]\subset S=R[x_n]$ be generated by
forms $\gg=\{g_0,\ldots,g_{n-1}\}$ of dgree $d\geq 1$.
Let
$$\cdots\rar\bigoplus_{j=1}^{m_1} R(-a_{1j})\xrightarrow{\phi} \bigoplus_{i=0}^{n} R(-d) \xrightarrow{\gg}R \rightarrow R/I\rightarrow 0$$
and
 $$\cdots\rar\bigoplus_{j=1}^{s_1} S(-C_{1j})\xrightarrow{\psi}\bigoplus_{j=1}^{s} S(-C_{j})\xrightarrow{\pi}S \rightarrow S/IS:f\rightarrow 0$$
 stand for minimal graded free resolutions of $R/I$ and $R/IS:f$, over $R$ and $S$ respectively,
from which we trivially derive minimal graded free resolutions of $S/qIS$ and $R/q(IS:f)$ over $S$:
$$\cdots \bigoplus_{j=1}^{m_1} S(-a_{1j}-\deg(q))\xrightarrow{\phi_1=\phi} S(-(d+\deg(q)))^{n}
\xrightarrow{q\,\gg}S \rightarrow S/qIS\rightarrow 0 ,$$
 $$\cdots\rar\bigoplus_{j=1}^{s_1} S(-C_{1j}-\deg(q))\xrightarrow{\psi_1=\psi}\bigoplus_{j=1}^{s}
 S(-C_{j}-\deg(q))\xrightarrow{q\,\pi}S \rightarrow S/q(IS:f)\rightarrow 0.$$
Shifting the second of these resolutions by $-(d+\deg(q))$, one obtains a map of complexes, where the vertical
 homomorphisms are also homogeneous of degree $0$
induced by multiplication by $f$ on the rightmost modules:
{\scriptsize
 $$
 \begin{array}{ccccccccccccc}
 \cdots\kern-6pt &\rar \kern-6pt & \bigoplus_{j=1}^{m_{i}} S(-a_{ij}-\mathfrak{d})\kern-6pt &\rar
 &\kern-6pt \cdots \kern-6pt & \stackrel{\phi}{\rar}\kern-6pt  &  S(-(d+\mathfrak{d}))^{n}
 \kern-6pt & \rar\kern-6pt & \kern-6pt S \kern-6pt & \rar   \kern-6pt& S/qIS  \kern-6pt &\rar  & 0 \\[5pt]
  && c_i(f)\uparrow &&&& c(f)\uparrow && \cdot f\uparrow && \cdot f\uparrow &&\\[5pt]
   \cdots\kern-6pt &\rar \kern-6pt & \bigoplus_{j=1}^{s_{i}} S(-C_{ij}-(d+2\mathfrak{d}))\kern-6pt &\rar
   &\kern-6pt \cdots \kern-6pt & \kern-6pt\stackrel{\psi}{\rar}  &
 \kern-6pt \bigoplus_{j=1}^{s} S(-C_{j}-(d+2\mathfrak{d})) &\kern-6pt \rar & \kern-6pt S(-(d+\mathfrak{d}))
 & \rar \kern-6pt  & \frac{S}{q(IS:f)}(-(d+\mathfrak{d})) & \kern-6pt\rar  & 0
 \end{array}
 $$
 }
 where we have written $\mathfrak{d}:=\deg(q)$ for editing purpose.

\begin{Proposition}\label{mapping_cone}
The mapping cone of the above map of complexes is a graded free resolution of the ideal $(qI,f)$
of $(q,f)$-type with underlying ideal $I$.
Moreover, if hd$(S/q(IS:f))\leq {\rm hd}(R/I)-1$
{\rm(}e.g., if $f\in IS$ and $I$ has codimension $\geq 2${\rm)}
then ${\rm hd}(S/(qI,f))\leq {\rm hd}(R/I).$
\end{Proposition}
\demo This result was essentially proved in \cite{HSstellar}.
The present case only requires minor changes.
\qed

\smallskip

Note the peculiarity of the resolution where $q$ intervenes by way of degree shift, while $f$ is travestied
in the generators of the quotient $IS:f$.
In particular, we draw attention to the syzygy matrix of the generators of $(qI,f)$, which has the form
\begin{equation}\label{syzygies}
\Psi= \left(
\begin{array}{cc}
\phi & c(f)\\
\mathbf{0} & -q\pi
\end{array}
\right).
\end{equation}
Here $\phi$ denotes a syzygy matrix of the given set of generators of $I$, while $\pi:S^s\surjects IS:(f)$ stands for
a surjective $S$-module homomorphism based on the given homogeneous generators $\{c_1,\ldots,c_s\}\subset S$  of $IS:(f)$,
and $c(f):S^s\lar S^{n+1}$ is the induced {\em content map} whose $j$th column vector gives the coefficients
of $fc_j$ as a combination over $S$ of the generators of $I$.

\smallskip

We will write $D:=\deg(f)=d+\mathfrak{d}$.

\begin{Example}\label{inclusion}\rm
Suppose that $f\in IS$. Then $\pi$ is the identity map of $S$ and the content map $c(f):S\rar S^{n+1}$ is represented
by one single column. A graded free resolution of $(qI,f)$ has the form
$$0\rar F_r\lar \cdots \lar F_2\lar F_1\oplus S(-(D+\mathfrak{d}))\lar S^{n+1}(-D)\lar S,$$
where $0\rar F_r\rar \cdots \rar F_2\rar F_1\rar S^n(-D)\rar S$ is a graded free resolution of $qIS$ over $S$
with suitable self-understood shifts.
\end{Example}
Two important cases emerge as follows.

\subsubsection{Hilbert--Burch ideal of $(q,f)$-type}

\smallskip

If $I$ is a codimension $2$ Cohen--Macaulay ideal then so is $(qI,f)$.
Namely,  a graded free resolution of $(qI,f)$ is
$$0\rar \left(\bigoplus_j S(-(\mathfrak{d} +d_j))\right)\oplus S(-(D+\mathfrak{d}))\lar S^{n+1}(-D)\lar S,$$
where $0\rar \bigoplus_j S(-d_j)\lar S(-d)^n\lar S$ is a graded free resolution of $IS$.

Unfortunately, the classification of codimension $2$ Cohen--Macaulay Cremona maps is unknown.
If $n-1=2$, i.e., for plane Cremona maps, this property is equivalent to requiring that the base ideal be
saturated. It has been proved in \cite[Theorem 1.5]{HS} that a plane Cremona map of codimension $2$ and degree
at most $4$ is Cohen--Macaulay. This is false for degree $\geq 5$ and the classification gets harder and
harder as the degree increases.

\smallskip

The Cremona map of degree $3$ in (\ref{the_cubic}) falls within this class: it is based on the standard plane quadratic
Cremona map, which is Cohen--Macaulay of codimension $2$.

\subsubsection{Almost Koszul ideal of $(q,f)$-type}

This is the case where $I$ is generated by a regular sequence (necessarily of degree $1$), which is resolved by the Koszul complex
on the elements of the regular sequence.
The graded free resolution of $(qI,f)$ has the form
$$0\rar S(-(n +\mathfrak{d}))\rar \cdots \rar S(-(3 +\mathfrak{d}))^{{n}\choose {3}}\rar S(-(2 +\mathfrak{d}))^{{n}\choose {2}}
\oplus S(-(D+\mathfrak{d}))\rar S^{n+1}(-D)\rar S.$$

This class of maps and the Hilbert--Burch ones are quite apart.
Indeed, the only Cremona map of de Jonqui\`eres which are both Cohen--Macaulay of codimension $2$ and
based on a complete intersection is a classical plane de Jonqui\`eres map (thus forcing $n=2$).

Infringing our notation for a minute, writing $d:=\deg(q)+1=\mathfrak{d}+1$,
the graded free resolution of an almost Koszul Cremona map (i.e., based on a complete intersection)
has the form
$$0\rar S(-(d +n-1))\rar \cdots \rar S(-(d +2))^{{n}\choose {3}}\rar S(-(d +1))^{{n}\choose {2}}
\oplus S(-(2d -1))\stackrel{\phi}{\lar} S^{n+1}(-d)\rar S.$$

\subsection{Inverse results}

From the previous subsection we transcribe the relevance of the two discussed special classes of ideals of
$(q,f)$-type in the following condensed result:

\begin{Proposition}
The de Jonqui\`eres group $\star_o(1;\PP^3)$ is generated by de Jonqui\`eres maps of Hilbert--Burch type
and of almost Koszul type.
\end{Proposition}
In this part we will reverse our considerations, by assuming a resolution format is given from which we
wish to deduce the nature of the ideal.
Unfortunately, due to the difficulty in classifying the plane Cremona maps with Cohen--Macaulay base ideals, it becomes
hard to recover which elements of $\star_o(1;\PP^3)$ are determined by the minimal free resolution of the
respective perfect base ideals.
In the subsequent part we will find a satisfactory result for the case of almost Koszul maps.

\subsubsection{Homological characterization}

Since we expect, under suitable circumstances, to recover a de Jonqui\`eres type of ideal, we will call the
base ring $S$ (instead of $R$) and denote the ideal by $J$ (instead of $I$).
We will, however, for light reading stick to $d$ (instead of $D$) for the degree of the generating forms.

Thus, let $S=k[x_0,\ldots,x_n]$ be a standard graded polynomial ring over a field $k$, where $n\geq 3$, and let $J\subset S$ denote
a homogeneous ideal of codimension $\geq 2$,  having a minimal free graded resolution of the form
\begin{equation}\label{resolution}
0 \rar S(-(d+2)) \stackrel{\psi}{\lar} S^n(-(d+1))\oplus S(-(2d-1)) \stackrel{\phi}{\lar} S^{n+1}(-d) \lar J \rar 0.
\end{equation}
Clearly, for such homological dimension, $J$ is saturated (i.e., $J^{\rm sat}:=J:(\xx)^{\infty}=J$)
and $d\geq 2$.
Moreover, if $n=3$ then $J$ has codimension $\leq 2$. In fact, otherwise $J$ would have codimension $\geq 3$, hence of
codimension exactly $3$ because
the projective dimension is $3$. Then $S/J$ would be Cohen--Macaulay, hence Gorenstein because the
Cohen--Macaulay type is $1$. But this is nonsense since any (homogeneous or local) Gorenstein ideal of codimension $c\geq 3$ is minimally
generated by at least $c+2\geq 5$ elements.

Thus, $J$ has codimension $2$ if $n=3$.
We assume throughout that $n=3$.
Recall that the {\em unmixed part} $J^{\rm un}$ of an ideal $J\subset S$ is the intersection of its primary components of minimal codimension and that $J$ is said to be {\em unmixed}
if $J^{\rm un}=J$.

\begin{Theorem}\label{main_homological}
Let $n=3$ and let $J\subset S$ denote an ideal of codimension $2$ having a minimal free resolution as in
{\rm (\ref{resolution})}.
If $d=2$ we assume, moreover, that $J$ is not unmixed.
Then
\begin{enumerate}
\item[{\rm (i)}]  $S/J$ has a unique associated prime of codimension $3$ and, moreover, this prime is generated
by the entries of $\psi\,${\rm ;}
\item[{\rm (ii)}] Up to a linear change of variables and elementary row operations one has
\end{enumerate}
\[\psi=\left(
\begin{array}{c}
x_0\\
x_1\\
x_2\\
0
\end{array}
\right),\quad
\phi=\left(
\begin{array}{ccc@{\quad\vrule\quad}c}
&&&\raise5pt\hbox{$-q_0$}\\
&&&\raise2pt\hbox{$-q_1$}\\
&\raise8pt\hbox{$\mathcal{K}$}&&\raise5pt\hbox{$-q_2$}\\
\multispan4\hrulefill\\[1pt]
0&0&0&q
\end{array}
\right),\quad
J=(qx_0,qx_1,qx_2,q_0x_0+q_1x_1+q_2x_2),
\]
for suitable forms $q_0,q_1,q_2,a\in S$ of degree $d-1$, where $\mathcal{K}$ is the Koszul syzygy
matrix of the regular sequence $x_0,x_1,x_2\,${\rm ;} in particular, $J^{\rm un}=(q,q_0x_0+q_1x_1+q_2x_2)$,
a complete intersection of degree $(d-1)d$.
\end{Theorem}
\demo
Let us first assume that $d\geq 3$.
Then $2d-1\geq d+2$, hence we may assume that the last entry in $\psi$ vanishes.
Since the remaining entries are linear and generate an ideal of codimension $\geq 3$ by
the acyclicity criterion of Buchsbaum--Eisenbud (\cite[Theorem 20.9]{E}), then they form
a regular sequence of linear forms. By a change of variables, we may assume that the entries of
$\psi$ are $x_0,x_1,x_2,0$.
Consider the $4\times 3$ linear submatrix $\mathcal{L}$ of $\phi$. Its rows are Koszul
relations of $x_0,x_1,x_2$, therefore it has rank $\leq 2$.
Clearly then its rank is exactly $2$ since $\phi$ has rank $3$.
Therefore, we may up to elementary row operations write

$$\phi=\left(
\begin{array}{ccc@{\quad\vrule\quad}c}
&&&\raise5pt\hbox{$-q_0$}\\
&&&\raise2pt\hbox{$-q_1$}\\
&\raise8pt\hbox{$\mathcal{K}$}&&\raise5pt\hbox{$-q_2$}\\
\multispan4\hrulefill\\[1pt]
0&0&0&q
\end{array}
\right),
$$
where $\mathcal{K}$ denotes the transposed Koszul syzygy matrix of $x_0,x_1,x_2$.
Note that all the operations so far have changed the entries on the right most column
of $\phi$, but not their degrees ($=d-1$).

But this immediately implies that the columns of $\mathcal{K}$ are syzygies of say, the first three minimal generators
of $J$. This implies that the latter have the form $px_0,px_1,px_2$, for some form $p\in S$ of degree $d-1$.
This already proves the assertion in (i), namely, since $p\notin J$ by degree consideration, the prime ideal
$P=(x_0,x_1,x_2)$ is an associated prime of $S/J$.
Moreover, by a well-known fact (see, e.g., \cite[Corollary 20.14(a)]{E}) this is necessarily the only
associated prime thereof of codimension $3$.

To complete the proof of (ii) it remains to show that we can assume $p=q$ and the fourth minimal generator of $J$ is of
the stated form.
Let $\mathcal{N}$ denote, say, the submatrix of $\phi$ consisting of the three rightmost columns, which we may clearly
assume has rank $3$.
Recall that, up to order and signs,
the set of $3\times 3$ minors of $\mathcal{N}$ divided by their $\gcd$ coincide with the given set of minimal generators
of $J$ on which (\ref{resolution}) is based.
This fact is well-known and follows by dualizing (\ref{resolution}) into $S$ to get
$$0\rar J^*\simeq S \lar S^4(d) \stackrel{\phi^t}{\lar} S^3((d+1))\oplus S((2d-1)),$$
from where follows that the entries of the vector generating the image of $S \lar S^4(d)$ are the maximal minors of the
transpose of $\mathcal{N}$ divided by their $\gcd$.

In the present case, by the explicitness of $\mathcal{N}$, we immediately see that this $\gcd$ is $x_0$, so the
determinant of
$$\left(
\begin{array}{ccc}
-x_0 & 0 & -x_2\\
0 & -x_0 & x_1\\
q_2 & q_1 & q_0
\end{array}
\right),
$$
further divided by $x_0$ gives the required expression up to signs adjustment.
To verify that $p$ can be taken to be $q$, one inspects easily the other $3$-minors.

The last statement in (ii) is obvious since localizing at any minimal prime of codimension $2$, a variable
among $x_0,x_1,x_2$ becomes invertible.

\medskip

We now consider the case $d=2$.
Our starter fails right at the outset since there are in fact unmixed ideals with the given resolution
shape -- e.g., the classical $J=(x_0,x_1)\cap (x_2,x_3)$.
Thus, we must assume that $J$ is not unmixed.
But then $S/J$ has some associated prime of codimension $3$.
As remarked earlier,  any such prime must contain the ideal generated by the
entries of $\psi$. Since these are all linear, we conclude as before that there is only one associated
prime of codimension $3$, necessarily generated by a regular sequence of $3$ linear forms.
This shows (i) for $d=2$ and we can pick up from here by repeating the argument used in the case where $d\geq 3$.
\qed

\medskip

We say that two rational maps  $\mathcal{F}, \mathcal{G}: \PP^n\tor\PP^n$ are {\em linearly equivalent} if
they belong to the same orbit
of the action of $\pgl(n+1,k)\times \pgl(n+1,k)$ (``source'' - ``target''), i. e., if there exist projective
transformations $\ell_1, \ell_2: \PP^n\tor\PP^n$ such that ˜$\mathcal{G} = \ell_1\circ\mathcal{F}\circ\ell_2$.
Clearly, the base ideals of two linearly equivalent maps have the same
 algebraic invariants -- most certainly are generated in the same degree, but moreover admit equivalent minimal graded
 free resolutions and the same corresponding graded Betti numbers.

As a consequence of the above results, we obtain a homological characterization of the subgroup $\star_o(n-1;\PP^3)$ of
Proposition~\ref{pro_fact} for $n=3$.

\begin{Theorem}
Let $\mathcal{F}$ denote a Cremona map of  $\PP^3$ of degree $d\geq 2$.
The following conditions are equivalent:
\begin{itemize}
\item[{\rm (a)}]  $\mathcal{F}$ is linearly equivalent to an element of $\star_o(2;\PP^3)\,${\rm ;}
\item[{\rm (b)}] The base ideal $J\subset S=k[x_0,x_1,x_2,x_3]$ of $\mathcal{F}$ admits a minimal
graded free resolution of the form
\[\xymatrix{0\ar@{->}[r]&S(-(d+2))\ar@{->}[r]&S^3(-(d+1))\oplus S(-(2d-1))\ar@{->}[r]&S^4(-d)\ar@{->}[r]&S/J},\]
and is not unmixed if $d=2$.
\end{itemize}
\end{Theorem}

We list some further properties of the above free resolution of the base ideal of a map $\mathfrak{F}$
in $\star_o(2;\PP^3)$.
By Theorem~\ref{main_homological}, we may assume that $J=(q(x_0,x_1,x_2),f)$, for some relatively prime forms $q,f\in S$
of respective degrees $d-1,d$, with $f\in (x_0,x_1,x_2)S$.
Moreover, $\mathfrak{F}$ being birational implies that $q$, and hence $f$ as well, is an $x_3$-monoid.
Let us write $f=\alpha x_3+\beta, q=\gamma x_3+\delta$,
with $\alpha,\beta,\gamma, \delta\in k[x_0,x_1,x_2]$ forms of degrees $d-1, d, d-2,d-1$, respectively.

\begin{Proposition}\label{details}
With the above notation, set $P=(x_0,x_1,x_2)$, Then:
\begin{enumerate}
\item[{\rm (1)}] The following conditions are equivalent:
\begin{enumerate}
\item[{\rm (i)}] $P$ is an embedded  prime of $S/J$
\item[{\rm (ii)}] Either $\mathfrak{F}$ has degree $\geq 3$  {\rm (}i.e., $\gamma\neq 0$ and $\deg (\gamma)\geq 1${\rm )}
or else $\gamma=0$.
\item[{\rm (iii)}] $(\alpha x_3+\beta, \gamma x_3+\delta)\subset P$.
\end{enumerate}
\item[{\rm (2)}] If the ideal $(c,\beta,\gamma, \delta)$ containing $(\alpha x_3+\beta, \gamma x_3+\delta)$ is a
proper ideal of codimension $\geq 3$ then the minimal primes of $S/(\alpha x_3+\beta, \gamma x_3+\delta)$ are
those that contract in $k[x_0,x_1,x_2]$ to an irreducible factor of the determinant $\alpha\delta-\beta\gamma\,${\rm ;}
otherwise {\rm (}i.e., if $(\alpha,\beta,\gamma, \delta)$ has codimension $2${\rm )}
the additional minimal primes are the defining primes of straight lines through $o$ and a point of
$V(\alpha,\beta,\gamma, \delta)\cap  V(x_3)$.
\end{enumerate}
\end{Proposition}
\demo
(1) Set $q:=\gamma x_3+\delta$.

(i) $\Rightarrow$ (ii) If $\mathfrak{F}$ has degree $2$ and $\alpha \neq 0$, then $q\notin P$. But if $P$ is embedded
then $P\supset Q$ for some minimal prime
of $S/J$ of codimension $2$. Since $q\notin Q$ then $Q$ must contain the variables $x_0,x_1,x_2$, hence contains $P$; a contradiction.

\smallskip

(ii) $\Rightarrow$ (iii) This is obvious.

\smallskip

(iii) $\Rightarrow$ (i) If $P$ is a minimal prime of $S/J$ and If $(\alpha x_3+\beta, \gamma x_3+\delta)\subset P$
then $P$ contains a minimal prime
of $S/(\alpha x_3+\beta, \gamma x_3+\delta)$; since the latter has codimension $2$ it is also a minimal prime of
$S/J$ because clearly $J\subset (\alpha x_3+\beta, \gamma x_3+\delta)$.
Therefore $P$ is an embedded prime of $S/J$.

\smallskip

(2) The ideal $(\alpha ,\beta,\gamma, \delta)$ is proper if and only if $\deg (\gamma)>0$.
Therefore, according to item (1)  this ideal is proper if and only if $P$ is an embedded prime of $S/J$.
If, moreover, $(\alpha ,\beta,\gamma, \delta)$ has codimension $\geq 3$ then the minimal primes of
$S/(\alpha ,\beta,\gamma, \delta)$ cannot be minimal primes of $S/(\alpha x_3+\beta, \gamma x_3+\delta)$.
Then, locally at any minimal prime $Q$ of $S/J$,  $\alpha \delta-\beta\gamma$ is a  generator of
$(\alpha x_3+\beta, \gamma x_3+\delta)$.
This implies that some irreducible factor of $\alpha \delta-\beta\gamma$ is a minimal generator of $Q$,
hence generates its contraction to $k[x_0,x_1,x_2]$
since this contraction has codimension $1$.
The alternative statement follows immediately since $\alpha ,\beta,\gamma, \delta\in k[x_0,x_1,x_2]$ now
implies that every minimal prime of $S/(\alpha ,\beta,\gamma, \delta)$
is generated by two linear forms in $k[x_0,x_1,x_2]$ defining the straight line through a point of
$V(\alpha ,\beta,\gamma, \delta)\cap  V(x_3)$ and $o$.
\qed

\subsubsection{Homologically near cases}

The purpose of this piece is to convey examples of ideals whose minimal free resolution is obtained from (\ref{resolution}) by a slight perturbation of its twists. Such ideals will actually be base ideals of Cremona maps which, therefore,
cannot be elements of $\star_o(2;\PP^3)$.

The first example appears on M. Noether's original papers (\cite{Noe71bis}; also
\cite[Example after Remark 2.3]{cremona}).

\begin{Example}
$J=(x_0x_3, x_1x_3, x_0(x_1-x_2),x_1(x_0-x_1))$.
\end{Example}
The minimal free resolution is of the form
$$0\rar S(-5)\lar S^3(-3)\oplus S(-4) \lar S^4(-2) \lar J\rar 0.$$
Note that it fits the template
$$0\rar S(-(d+3))\lar S^3(-(d+1))\oplus S(-2d) \stackrel{\phi}{\lar} S^4(-d) \lar J\rar 0.$$
Here  the linear submatrix of $\phi$ has rank $3$ (not $2$ as in the de Jonqui\`eres case) and the coordinates of  the tail map
generate a radical ideal in degrees $1,2$ whose minimal primes are
$$(x_0,x_1,x_3),\, (x_1,x_2,x_3),\, (x_0-x_2,x_1-x_2,x_3).$$
The unmixed radical of $J$ is $(x_0, x_1)$, hence $J$ has one embedded associated prime and two minimal primes of codimension $3$.

Moreover, a calculation with {\em Macaulay} (\cite{Macaulay}) shows that the initial degree of $J:P$ is $2$,
where $P$ is any of the above three associated primes. Therefore, there is no associated prime of codimension $3$
driven inside $J$ by a form of degree $1$.

\medskip

The next example defines the polar map of the determinant of a so-called $3\times 3$ sub-Hankel matrix
(\cite{CRS}, also \cite{MarAron}):
$$\phi=\left(
\begin{array}{ccc}
x_0 & x_1 & x_2\\
x_1 & x_2 & x_3\\
x_2 & x_3 & 0
\end{array}
\right).
$$

\begin{Example}
$J=(x_3^2, x_2x_3, x_2^2-2/3x_1x_3, x_1x_2-x_0x_3)$.
\end{Example}
The minimal free resolution is again of the form
$$0\rar S(-5)\lar S^3(-3)\oplus S(-4) \lar S^4(-2) \lar J\rar 0.$$
This time around, however, $S/J$ admits a unique associated prime of codimension $3$ and this prime
is an embedded prime.
This makes up for a sensitive geometric distinction between this example and the previous one
(although both have the same degree ($=1$) as schemes): the first is the scheme-theoretic union
of a straight line with an embedded point and two isolated points, while the second is a straight
line with an embedded point.
Clearly, any of these two is very distinct from de Jonqui\`eres ideal, whose scheme has
multiplicity $2$. Nevertheless, the second example is sort of more akin to a de Jonqui\`eres
as it is in a sense an ``iteration'' of Cohen--Macaulay de Jonqui\`eres schemes \cite[Remark 4.6 (b)]{CRS}),

\noindent {\bf Authors' addresses:}

\medskip

\noindent {\sc Ivan Pan}, Centro de Matem\'atica, Facultad de Ciencias, Universidad de la Rep\'ublica\\
 11400  Montevideo, Uruguay\\
{\em e-mail}: ivan@cmat.edu.uy\\

\smallskip

\noindent {\sc Aron Simis},  Departamento de Matem\'atica, CCEN, Universidade Federal
de Pernambuco\\
 50740-560 Recife, PE, Brazil\\
 and\\
 Departamento de Matem\'atica, CCEN, Universidade Federal
da Para\ii ba\\
58059-900 Jo\~ao Pessoa, PB, Brazil.\\
{\em e-mail}:  aron@dmat.ufpe.br

\end{document}